\documentclass[12pt]{article}
\usepackage{graphicx,amssymb, amsmath, vmargin}
\setpapersize{custom}{21cm}{29.7cm}
\setmarginsrb{2cm}{1cm}{2cm}{3.5cm}{0pt}{0pt}{0pt}{0pt}
\author{Mathieu Dutour}
\title{Computational methods for cones and polytopes with symmetry}

\newtheorem{theor}{Theorem}

\begin{document}
\newcommand{\R}{\ensuremath{\mathbb{R}}}
\newcommand{\N}{\ensuremath{\mathbb{N}}}
\newcommand{\Q}{\ensuremath{\mathbb{Q}}}
\newcommand{\C}{\ensuremath{\mathbb{C}}}
\newcommand{\Z}{\ensuremath{\mathbb{Z}}}
\newcommand{\T}{\ensuremath{\mathbb{T}}}

\maketitle

\begin{abstract}
Every polyhedral cone can be described either by its facets or by its 
extreme rays. Computation of one description from the other is 
a problem that can be very complex, i.e. one encounter the combinatorial 
explosion. We present here several methods which allow us to use symmetries
for the computation of extreme rays and facets of cones.
I give an algorithm for computing adjacencies using a linear programming 
approach. I give a complete description of the method of incidence, adjacency decomposition method and of the subcone method used for computing with symmetric cones.
\end{abstract}

\section{Introduction and motivation}
\noindent A {\it cone} in $\R^n$ is a set stable by addition and multiciplication by $\lambda\in\R_+$. A {\it polyhedral cone} is a cone defined by a finite number of linear inequalities. It is a theorem of Minkovski that those polyhedral cones can also be described as positive span of a finite number of {\it rays}. A ray is called an {\it extreme ray} if it is not the sum of two different rays. An inequality is called a {\it facet} if it cannot be expressed as a sum of two different inequalities.

So, there is essentially two description for a cone, either by extreme rays or by facets. The main problem of this paper is given a description by facets to obtain a description by extreme rays. The dual problem of computing a description by facets from a description by extreme rays is equivalent to the previous one.

In many mathematical, physical, and practical problems we have a description by facets (or extreme rays), the interesting objects being the extreme rays (or facets). So, the problem of translating the descriptions is a very interesting one and many programs exist that do it (cdd (\cite{Fu}), porta \footnote{http://www.zib.de/Optimization/Software/Porta}, lrs (\cite{Av}), etc).

But, in many cases one encountered the so called, combinatorial explosion, i.e. while the description by facets (or extreme rays) remain simple, the dual description has a huge number of extreme rays (or facets).

These cones have in many case some symmetry groups. Our goal in this paper is to present several methods for computing extreme rays and facets of cones using their symmetry group.

Two classical algorithms exist for computing extreme rays:
\begin{enumerate}
\item The Double description method, also called Fourier Motzkin method, consists of taking an initial list of $n$ facets of the cone and computing extreme rays of it (it is a simplicial cone so this is immediate), then to add the remaining facets one by one in order to have the list of extreme rays. The method is called Double Description since at any time of the process we know the description by facets and by extreme rays. This method is the one of cdd by Komei Fukuda and Avis (see \cite{FM} and \cite{FP}). This method has some variants (see for example \cite{VB})
\item The Pivoting method consists of taking an extreme ray of the cone and finding extreme rays adjacent to it. This method is generally inefficient for extreme rays with high degeneracy that occur in combinatorial optimization.
\end{enumerate}

The cones and polytopes appear in many context like Metric Geometry (metric cone, quasi-metric cone, super-metric cone (see \cite{DL}, \cite{DD-dm}, \cite{La}, \cite{GrMet}), Max-cut problem (cut cone and cut polytope, see \cite{DL} and \cite{D-L}), Geometry of numbers (Voronoi polytope and Voronoi Domain, see \cite{Jaq}), Combinatorial optimization (Traveling Salesman Polytope and Linear Ordering Polytope, see \cite{CR}), Game theory (solitaire cones, see \cite{AD}), Physics (see \cite{BCLM} and \cite{Er}).

%The cones appear in several contexts, we list some of them
%\begin{enumerate}
%\item The metric cones, quasi-metric cones, supermetric cones are described by their facets (see \cite{DL}, \cite{DD-dm}, \cite{La}, \cite{GrMet}). There also exists a dual theory of cut cones $CUT_n$ and some variants (see \cite{DDP} and \cite{DL}). All these cones have symmetry group $Sym(n)$.

%\item The Voronoi Cell of a lattice is described by its facets, its vertices are the center of the Delaunay cells of the lattice (called holes in \cite{CS}). Another case is the Voronoi domain defined by extreme rays. The facets of this cone correspond to the adjacent perfect forms (see \cite{Jaq}).

%\item The hypermetric cone is a cone defined by the hypermetric inequalities, the extreme rays correspond to extreme Delaunay polytopes (see \cite{DD-hyp} and \cite{DL}).

%\end{enumerate}
There are two main kinds of symmetry groups: Permutation groups like $Sym(n)$ arising in combinatorial optimization and matrix groups arising in geometry of numbers (see \cite{Jaq}).

Note also that the problems of finding the size of an orbit and of finding canonical representatives may be very difficult (see for example \cite{LN}, \cite{Jaq}). Althought we use the simpler language of cones here, the methods apply also to polytopes.

\section{General infos on cones and polynomially solvable questions}
\noindent A cone in $\R^n$ can be defined in several ways, first as an intersection of a finite number of half spaces:
\begin{equation*}
x\in\R^n\mbox{~such~that~}\sum_{i=1}^{n}a_{i,k}x_i\geq 0,\,\,\forall 1\leq k\leq m
\end{equation*}
or as a set spanned by vectors:
\begin{equation*}
x\in\R^n\mbox{~such~that~}\exists \lambda_l\geq 0\mbox{~with~}x=\sum_{l=1}^{p}\lambda_l v_l\;.
\end{equation*}
We assume that the cones have full dimension $n$. 
Given $v \in \mathbb{R}^{n}$, the inequality $\sum_{i=1}^nv_ix_i \geq 0$ is 
said to be {\it valid} for $C$, if it holds for all $x \in C$. Then the set 
$\{ x \in C \vert \sum_{i=1}^nv_ix_i = 0 \}$ is called the {\it face of $C$, 
induced by the valid inequality $\sum_{i=1}^nv_i x_i \geq 0$}. A face of 
dimension $n - 1$ is called a {\it facet} of $C$; a face of dimension 
$1$ is called an {\it extreme ray} of $C$.

Two extreme rays of ${\cal C}$ are said to be {\it adjacent}, if they
span a two-dimensional face of ${\cal C}$. Two facets of ${\cal C}$ are 
said to be {\it adjacent}, if their intersection has dimension $n - 2$.
The {\it $1$-skeleton} graph of ${\cal C}$ is the graph 
whose nodes are the extreme rays of ${\cal C}$ and whose edges
are the pairs of adjacent extreme rays. 
The {\it$n-1$-skeleton} graph of ${\cal C}$ is the graph with node set being
the set of facets of ${\cal C}$ with an edge between two facets if
they are adjacent in ${\cal C}$.

Some face conditions may be useless to the cone description by face, they are called redundant. Some vector may also be redundant. In fact, we can find a finite subset of non-redundant faces of the cone ${\cal C}$ called {\em facets}. The non-redundant rays are called {\em extreme rays}.

Given a cone ${\cal C}$ we can define ${\cal C^*}$ as the cone whose extreme rays are facets vector of ${\cal C}$, the facets of ${\cal C^*}$ will be the extreme rays of ${\cal C}$. So, every method or theorem given here as a dual version with extreme rays and facets exchanged.

The set of faces of a cone ${\cal C}$ is a lattice denoted by $L({\cal C})$ under the inclusion. Given a set ${\cal F}$ of faces a face named $\sup\, {\cal F}$ is defined as the smallest face containing the elements of ${\cal F}$, while $\inf\,{\cal F}$ is defined as the biggest face contained in every element of ${\cal F}$.

The $k$ skeleton of ${\cal C}$ is the set of $k$ dimensional faces of ${\cal C}$. Two $k$-faces $F_1$, $F_2$ are adjacent if $\sup\,\{F_1, F_2\}$ has dimension $k+1$ and $\inf\,\{F_1, F_2\}$ has dimension $k-1$. The {\em $k$-skeleton} graph is the graph with node set being the set of $k$-faces of ${\cal C}$ with an edge between two $k$-faces if they are adjacent in ${\cal C}$.

{\it The linear programming way of testing adjacencies}
\begin{enumerate}
\item First take a set ${\cal F}=\{f_i\}_{i\in I}$ of inequalities $f_i(x)\geq 0$ and another inequality $f$, the question is: Is that $f$ redundant? reply is to consider the linear programming problem:
\begin{center}
minimize $f(x)$ in the cone defined by ${\cal F}$
\end{center}
if the minimum is $0$ then $f$ is redundant (i.e. it can be expressed as a sum with non-negative coefficients of elements of ${\cal F}$). If the minimum is $-\infty$ then $f$ is is non-redundant.
\item {\it Claim} Given a non-redundant set ${\cal F}$ of facets of a cone ${\cal C}$ and $(f,g)$ a pair of facets then $(f,g)$ are non-adjacent if and only if there exists $\lambda_i \geq 0$ such that
\begin{equation*}
\sum_{i\in {\cal F}-\{f,g\}} \lambda_i f_i=\lambda_f f+\lambda_g g
\end{equation*}
with $\lambda_f>0$, $\lambda_g>0$.
\item Using the above claim we have $(f,g)$ non-adjacent if and only if 
\begin{equation*}
\begin{array}{c}
\inf_{x\in {\cal C}} f(x)=0\\
\mbox{with~}{\cal C}=\{x\mbox{~such~that~}u(x)\geq 0,~u\in ({\cal F}-\{f,g\})\cup (-g)\}
\end{array}
\end{equation*}
\end{enumerate}

This method allow us to compute adjacencies of facets of a cone defined by facets without knowing its extreme rays. It is so very efficient.

{\it The linear programming method for computing $k$-faces}
\begin{enumerate}
%\item {\it We assume we know the set of facets and of extreme rays} Here is a way of computing $\sup\,{\cal F}$ and $\inf\,{\cal F}$. Suppose that ${\cal F}=(F_i)_{i\in I}$. For each face $F_i$ determine the set $E_i$ of extreme rays incident to $F_i$.\\
%Now put ${\cal S}=\cup_{i\in I} E_i$ and ${\cal I}=\cap_{i\in I} E_i$.\\
%Now determine the set of all facets adjacent to ${\cal S}$ and ${\cal I}$, this gives the facet description of $\sup\,{\cal F}$ and $\inf\,{\cal F}$. The extreme rays description can be obtained as well.
\item We assume here that a cone is defined by a set ${\cal ER}$ of extreme rays.
\item The linear programming method allow us to compute edges from the set of extreme rays, i.e. we find from the set of extreme rays the set of edges. From this set of edge one would like to compute the set of adjacent edges, i.e. we want to compute the set of edges that span a $3$ face of the cone ${\cal C}$.\\
More generally generally the computation of adjacencies of $k$-faces amount to the computation of elements of the $k+1$-skeleton.
\item Note that the edges are $1$-simplex but the element of the elements of the $3$ skeleton can be nonsimplicial. This fact can be seen in several ways:
\begin{enumerate}
\item The $3$ skeleton of icosahedron is made of triangles while for dodecaedron it is made of pentagon.
\item $n=3$ is the first dimension where cones can have an arbitrary high number of extreme rays.
\item In a list of $3$ rays that are linearly dependent, at least one of them is redundant, while there exist list of $4$ rays that are linearly dependent but all irredundant.
\end{enumerate}
\item {\it Linear Programming way of testing the adjacency}: Let $F_1$ and $F_2$ be two $k$-faces, we assume that $\inf\{F_1, F_2\}$ has dimension $k-1$ (i.e. it is one of the $k-1$ faces found at the last step). The process is as follow:
\begin{enumerate}
\item We determine the set ${\cal E}$ of extreme rays $e$ such that $e\in Vect\,(F_1\cup, F_2)$
\item $\sup\{F_1\cup F_2\}$ has dimension superior to $k+1$ if and only if there exist $\lambda_{e'}\geq 0$ and $\lambda_e>0$ such that
\begin{equation*}
\sum_{e'\in {\cal ER}-{\cal E}} \lambda_{e'} e'=\sum_{e\in {\cal E}}\lambda_e e
\end{equation*}
this correspond to the linear programming
\begin{equation*}
\begin{array}{c}
\inf_{x\in {\cal C}} es.x=0\\
\mbox{with~}{\cal C}=\{x\mbox{~such~that~}e.x\geq 0\,\forall e\in {\cal ER}-{\cal E}\mbox{~and~}e.x\leq 0\,\forall e\in {\cal E}-\{es\}\}
\end{array}
\end{equation*}

\end{enumerate}

\end{enumerate}

\section{The method of Symmetry}
\noindent Consider a cone ${\cal C}$ with a linear symmetry group $G$ acting on it. If $F$ is a facet then the symmetry group of the facet is
\begin{equation*}
Sym_F=\{g\in G\mbox{~such~that~}g(F)=F\}
\end{equation*}
is a subgroup of $G$ satisfying to $|Sym_F|\times |O_F|=|G|$ with $O_F$ the orbit of $F$. Furthermore if $F$ and $F'$ belong to the same orbit then $Sym_F$ and $Sym_{F'}$ are conjuguate subgroups of $G$. The conjugacy classes of $Sym(n)$ are
\begin{equation*}
\begin{array}{|c|c|c|c|c|c|c|c|c|c|}
\hline
n&3&4&5&6 &7 &8 &9 &10&11\\
\hline
|\mbox{Conjugacy~classes}|&3&5&7&11&15&22&30&42&56\\
\hline
\end{array}
\end{equation*}
If $H$ is a subgroup of $G$ then we define ${\cal C}_H=\{d\in{\cal C}, \forall h\in H, h(d)=d\}$. The cone ${\cal C}_H$ has a dimension smaller than $n$. If ${\cal C}$ is defined by extreme rays then the extreme rays of ${\cal C}_H$ can be found easily.

\begin{theor}
If $F$ is a facet of ${\cal C}$ having $H$ as a symmetry group then $F$ is a facet of ${\cal C}_H$.
\end{theor}
We set $U_F=\{d\in {\cal C}\mbox{~s.t.~}F(d)=0\}$.
If $F$ has $H$ symmetry, then $U_F$ has also $H$ as a symmetry group.\\
Now define for any ray of ${\cal C}$
\begin{equation*}
SU(d)=\sum_{h\in H} h(d)
\end{equation*}
Denote $SU(X)=\{SU(d)\mbox{~s.t.~}d\in X\}$ for any $X\subset OMCUT_n$. 
Then $SU(U_F)$ is a set of symmetric elements, which are incident to $F$. 
Moreover, since ${\cal C}_H=SU({\cal C})$, we have $SU(U_F)\subset {\cal C}$. 
We denote $n$ the dimension of ${\cal C}$ and $m$ the dimension of 
${\cal C}_H$. By hypothesis $F$ is a facet, so $U_F$ has dimension 
$n-1$. The mapping $d\rightarrow SU(d)$ decrease dimension by at 
most $n-m$; so, we get that $SU(U_F)$ has dimension $m-1$, i.e. 
$F$ is a facet of ${\cal C}_H$.
\\

Nevertheless the facets of the cone ${\cal C}_H$ are not necessarily facets 
of ${\cal C}$, i.e. after computing facets of ${\cal C}_H$, we need
to check if they are indeed facets of ${\cal C}$.

Using this method we can found the facets of a given cone having a given 
symmetry. Since the number of conjugacy classes of a group is not large
we can compute orbits having a given non-trivial symmetry.

The cone $CUT_7$ has the property that all its facets have a non-trivial 
symmetry group but this is a particular case and in general we can't expect 
to find all facets in this way.

\section{The method of incidence}
Let ${\cal C}$ be a cone of dimension $n$ with a symmetry group $G$ and a set ${\cal F}$ of facets partitionned into orbits $\{O_1, \dots, O_p\}$ with representatives $\{F_1, \dots, F_2\}$

If $e$ is an extreme ray then it is incident to at least one facet of ${\cal C}$. By the group action, we can assume it is incident to a representative $F_i$ of orbit $O_i$. This incidence yield a linear equality and this equality allow us to find a cone of dimension $n-1$ of which $e$ is an extreme ray. This idea can be extended to faces of dimension $n-k$ with $k$ small. Each such faces is a cone of dimension $n-k$.

This idea was used by Grishukhin in \cite{GrMet} on $MET_7$ and by Christof and Reinelt in \cite{CR} and by \cite{DMP} for the metric cone $MET_8$. The method is the following: We find a set of orbits of $n-k$ faces of ${\cal C}$ using the linear programming method of the preceding section. Then for each of these $n-k$-faces we find their extreme rays using classical program like cdd (\cite{Fu}).

One problem of the method is that a given extreme ray belong to many $n-k$-faces so the orbit of extreme rays found will belong to different orbit of $n-k$-faces.

\section{Adjacency decomposition method}
Take a cone ${\cal C}$ of dimension $n$ defined by a set ${\cal E}$ of 
extreme rays of ${\cal C}$. Now consider the set ${\cal E}_F$ of extreme 
rays incident to a given facet $F$. The facet $F$ can be considered as 
a cone ${\cal C}_F$ of dimension $n-1$ generated by the set ${\cal E}_F$ of 
extreme rays.

Furthermore the mapping
\begin{equation*}
\begin{array}{rcl}
(L(C), \subset)&\rightarrow&(L(C_F), \subset)\\
G&\mapsto&G\cap F
\end{array}
\end{equation*}
is a lattice morphism. In particular the facets of the cone ${\cal C}_F$ 
correspond to the ridges included in $F$. It is computationally easy to 
map a facet of ${\cal C}_F$ to a ridge $R$ of ${\cal C}$, to find the 
facet $G$ of ${\cal C}$ such that $R=G\cap F$ and so, to build facets 
of ${\cal C}$. 

So, we can from a facet of ${\cal C}$ find the adjacent facets to $F$, 
the complexity of the computation being the one of finding facets of 
the cone ${\cal C}_F$ defined by its extreme rays, this complexity is 
related to the incidence number of $F$. If it is equal to $n-1$ then 
the facet is simplicial and this is trivial, if the facet has incidence 
larger than $n-1$ then this becomes more difficult and we encounter again 
the combinatorial explosion while on a lesser form.

The algorithm is so:
\begin{enumerate}
\item Take a list of orbits of facets of the cone ${\cal C}_F$, find their incidence numbers and orbit size.
\item Take a facet from the list which is non-treated and find adjacent facets to it as indicated above
\item Extract from this list of facets the orbits and update the  list of orbits.
\item Go to step 2 if there remain some non-treated facets.
\end{enumerate}

This method was applied by Deza, Fukuda, Pasechnik, and Sato in \cite{DFPS} for the metric polytope $MET_8$ and $MET_9$, by Christof and Reinelt, in \cite{CR} for $CUT_8$, $CUT_9$, the linear ordering polytope $P_{LO}^8$ and the Traveling Salesman Problem $P_{TSP}^{8}$ and $P_{TSP}^{8}$, by Deza, Dutour and Pantaleeva in \cite{DDP} for $OMCUT_5$, $OMCUT_6$ and $QMET_6$, by Deza and Dutour in \cite{DD-dm} for supermetric cones $SMET^{m, s}_n$.

The result of this method fall in three cases:
\begin{enumerate}
\item The method find all orbits of facets, this is the case of $OMCUT_5$, $P_{TSP}^8$ for example
\item The method find a set of orbits of facets which is conjectured to be complete, this is the case of $HCUT^2_6$, $HCUT^4_7$, $CUT_8$, $MET_8$, $SMET^{2, 2}_6$
\item The method give only an estimate of the number of orbits, this is the case of $MET_9$, $CUT_9$, $P_{TSP}^9$, $OMCUT_6$, $SCUT^{2, 2}_6$.
\end{enumerate}

So, this method allow us to draw a frontier of the combinatorial explosion. In the second case, we find facets adjacent to orbits except for a few ones that have a large number of incident extreme rays, example
\begin{enumerate}
\item For $CUT_8$, the triangle facet, the pentagonal facet and the seven gonal facets have a too large incidence number (see \cite{CR})
\item For $HCUT^4_7$ and $HCUT^2_6$, the non-negativity facets and the simplex facets.
\item For $MET_8$ the cuts (see \cite{DFPS}).
\item For the Traveling Salesman Polytope, the Linear Ordering Polytope see \cite{CR}
\end{enumerate}
The complexity of these facets is usually extremely high, but one may want to prove that the found description are complete.

The mathematical object that appear is the quotient graph of the skeleton graph with elements, the orbits $O_i$ of facets. Two orbits $O_i$ and $O_j$ being adjacent if there is a facet $F_i\in O_i$ and a facet $F_j\in O_j$ such that $F_i$ and $F_j$ are adjacent facet in the cone ${\cal C}$.

The mathematical conjecture is that the quotient graph with those orbits removed is still connected. As far as we know no result of this kind as ever been proved. The only result available is the Balinski theorem which asserts that the skeleton graph of a $d$-polytope is $d$ connected, i.e. the removing of $d$ elements leave it connected. This result could not be applied in our cases since the size of the problematic orbits is larger than the dimension of the underlying space.

Another approach is possible: These orbits of facets have a very simple expression and a large symmetry group (for example the symmetry group of the triangle facet of $CUT_n$ is $Z_2\times Sym(n-3)$) so, the preceding approach can be applied recursively to the cone ${\cal C}_F$. This approach succeeds with $HCUT^4_7$, we intend to apply it to the cut cone $CUT_8$.

This approach is the one of Jaquet in \cite{Jaq}. [NEED TO CHECK THIS, THE JAQUET METHOD SEEMS QUITE INVOLVED AND SMART, HE REWRITES EVERYTHING]

In order to run the adjacency decomposition method, one needs to have an initial set of facets, among the possible method are:
\begin{enumerate}
\item Use of linear programming to generate vertices
\item Special kind of facets known by the theory of particular cones (zero-extension for metric cones for examples)
\item Partial cdd output obtained by running cdd or other program and 
\end{enumerate}

\section{The Sub-Cone method}
We give here a method that was applied to only one cone, i.e. the hypermetric cone on seven vertices $HYP_7$. The cone $HYP_7$ is defined by $3773$ facets called {\it hypermetric facets} which form $14$ orbits under the group $Sym(7)$ (see \cite{B70} and \cite{Ba})
\begin{enumerate}
\item We have computed the facets of the cone $CUT_7$. There was $36$ orbits of facets \cite{Gcut}. $10$ of them are hypermetric.
\item The cone $HYP_7$ is decribed by $14$ orbits of hypermetric facets, $10$ already appear in the list of $36$ orbits of facets
\item We define $26$ subcones
\begin{equation}
H_p=\{x\in HYP_7\mbox{~and~}p(x)\leq 0\}
\end{equation}
for every representative $p$ of the $26$ orbits of non-hypermetric facets. We remove the redundant facets by linear programming (see \cite{Fu})
\item The non-redundant description of $H_p$ are simplicial cones which were easily computed, so giving us description of extreme rays of $HYP_7$.
\end{enumerate}
The extreme rays of $HYP_7$ belong to two classes: the cuts coming from $CUT_7$ and the $26$ found extreme rays. Theses extreme rays are adjacent only to cuts, so the adjacency decomposition method is not useful here: The skeleton graph of $HYP_7$ without the cuts is totally disconnected.

The success of the method rely on several factor: the knowledge of facets of $CUT_7$ which is not something easy \cite{Gcut} and the fact that the cones $H_p$ are disjoint making the computing process reasonable. The cuts form a complete graph in $HYP_7$ but this fact solely does not explain why the $H_p$ are disjoint since $MET_7$ which is in the same situation as $HYP_7$ does not share this property. All attempt to use this method in other contexts were inefficient.


\begin{thebibliography}{99}

\bibitem[Av]{Av}
D. Avis, {\em A C implementation of the reverse search vertex enumeration algorithm}, School of Computer Science, McGill University, Montreal, Canada, 1993, program lrs*.c available from ftp://mutt.cs.mcgill.ca/pub/C/

\bibitem[AD]{AD}
D. Avis and A. Deza, {\em On the binary solitaire cone}, Discrete Applied Mathematics 115-1 (2001) 3-14.

\bibitem[B70]{B70}
E.P. Baranovski, {\em Simplexes of $L$-subdivisions of euclidean spaces}, Mathematical Notes, {\bf 10} (1971) 827-834.

\bibitem[BCLM]{BCLM}
L. B\'enin, C. Cummin, L. Lapointe and P. Mathieu, {\em Fusion bases as facets of polytopes}, http://www.arxiv.org/abs/hep-th/0108213

\bibitem[Ba99]{Ba}
E.P. Baranovskii, {\em The conditions for a simplex of $6$-dimensional lattice to be $L$-simplex}, (in russian) Nauchnyie Trudi Ivanovo state university. Mathematica {\bf 2} (1999) 18--24.

\bibitem[ChRe96]{CR}
T.Christof and G.Reinelt, {\em Combinatorial optimization and small polytopes}, Top (Spanish Statistical and Operations Research Society) {\bf 4} (1996)
 1--64.

%\bibitem[CS]{CS}
%J.H. Conway and N.J.A. Sloane, {\em Sphere Packings, Lattices and Groups}, volume 290 of {\em Grundlehren der mathematischen Wissenschaften}, Springer Verlag.

\bibitem[DeDe94]{DD}
A.Deza and M.Deza, {\em The ridge graph of the metric polytope and some 
relatives}, in T.Bisztriczky, P.McMullen, R.Schneider and A.Ivic Weiss eds.
 Polytopes: Abstract, Convex and Computational (1994) 359--372.

\bibitem[DDP]{DDP}
M.Deza, M.Dutour and E.I.Panteleeva, {\em Small cones of oriented semi-metrics}, in preparation.

\bibitem[DDFu96]{DDF}
A.Deza, M.Deza and K.Fukuda, {\em On Skeletons, Diameters and Volumes of 
Metric Polyhedra}, in Combinatorics and Computer Science, Vol. {\bf 1120}
 of Lecture Notes in Computer Science, Springer--Verlag, Berlin (1996)
 112--128.

\bibitem[DeLa97]{DL}
M.Deza and M.Laurent, {\em Geometry of cuts and metrics}, Springer--Verlag, Berlin 1997.

\bibitem[D-L]{D-L}
M.Deza and M.Laurent, {\em Applications of cut polyhedra}, Journal of computational and applied mathematics {\bf 55} (1994) 191-216, 217-247.

\bibitem[DD01]{DD-hyp}
M.Dutour and M.Deza, {\em The hypermetric cone on seven vertices},
submitted (2001) http://il.arXiv.org/abs/math.MG/0108177

\bibitem[DD01]{DD-dm}
M.Deza and M.Dutour, {\em Data Mining for cones of metrics, quasi-metrics, $m$-hemi-metrics and $(m,s)$-super-metrics with large number of extreme rays or facets}, submitted (2002) http://www.arxiv.org/abs/math.MG/0201011

\bibitem[DFPS]{DFPS}
A. Deza, K. Fukuda, D. Pasechnik, and M. Sato, {\em On the skeleton of the metric polytope},  Lecture Notes in Computer Science {\bf 2098} Springer-Verlag, Berlin (2001) 125-136.

\bibitem[DMP]{DMP}
A. Deza, T. Mizutani, and D. Pasechnik, {\em Enumerating faces of polyhedra with large symmetry group}, To appear in Optimization: Modeling and Algorithms 15. The Inst. Stat. Math. pp. 179-185.

\bibitem[Er]{Er}
R.M. Erdahl, {\em Representability conditions}, in R.M. Erdahl and V.H.Smith Eds., Density Matrices and Density Functionnals (Reidel, Dordrecht, 1987) 51-75.

\bibitem[LN]{LN}
F. L\"ubeck and M.Neunh\"offer, {\em Enumerating Large Orbits and Direct Condensation}, Experimental Mathematics, {\bf 10-2} (2001), 197-205.

\bibitem[La]{La}
M. Laurent, {\em Graphic vertices of the metric polytope. Graph theory and combinatorics} Discrete Mathematics, {\bf 151} (1996), 131-153.

\bibitem[VB]{VB}
E.Viterbo and E.Biglieri, {\em Computing the Voronoi Cell of a Lattice: The Diamond-Cutting Algorithm}, IEEE Transactions on informations theory, {\bf 42-1} (1996) 161-171

\bibitem[Fu95]{Fu}
K.Fukuda, {\em http://www.ifor.math.ethz.ch/\texttt{\~}fukuda/cdd\_home/cdd.html}, 

\bibitem[Gr92]{GrMet}
V.P.Grishukhin, {\em Computing extreme rays of the metric cone for seven
 points}, European Journal of Combinatorics {\bf 13} (1992) 153--165.

\bibitem[Gr90]{Gcut}
V.P.Grishukhin, {\em All facets of the cut cone $C_n$ for $n=7$ are known}, European Journal of Combinatorics {\bf 11} (1990) 115--117.

\bibitem[Jaq]{Jaq}
Jaquet, {\em \'Enum\'eration compl\`ete des classes de formes parfaites en dimension $7$}, Ann Inst Fourier 43 (1993), 21-55 et th\`ese de doctorat Neuchatel

\bibitem[FM]{FM}
T.S.Motzkin, H.Raiffa, G.L.Thompson, and R.M.Thrall, {\em The double description method}, In H.W.Kuhn A.W.Tucker, editors, Contributions to the theory of games, volume 2. Princeton University Press, Princeton, N.J., 1953

\bibitem[FP]{FP}
K.Fukuda and A.Prodon, {\em Double description method revisited}, Combinatorics and computer science (Brest 1995) 91-111, Lecture notes in computer science 1120, Springer 1996


\end{thebibliography}
\end{document}